\documentclass[10pt,a4paper,twoside,leqno]{article}
\usepackage{amsfonts,amssymb}
\usepackage{eufrak}
\usepackage{amscd}

\newtheorem{defn}{Definition}

\vspace{0.5cm}

 4
\font\ebf=cmbx8
\font\erm=cmr8

\usepackage{graphicx}

\begin{document}

\thispagestyle{empty}

\noindent {\bf COMBINATORIAL INTERPRETATION\\ OF FIBONOMIAL
COEFFICIENTS }

\vspace{0.7cm} {\it Andrzej K. Kwa\'sniewski}

\vspace{0.7cm}

{\erm Higher School of Mathematics and Applied Informatics}

\noindent{\erm  Kamienna 17, PL-15-021 Bia\l ystok , Poland}

\vspace{0.5cm}

\noindent {\ebf Summary}

{\small   A  classical-like combinatorial interpretation of the
Fibonomial coefficients is provided following  [1,2]. An adequate
combinatorial interpretation of recurrence satisfied by Fibonomial
coefficients is also proposed. It is considered to be - in the
spirit classical- combinatorial interpretation like binomial
Newton and Gauss $q$-binomial coefficients or Stirling number of
both kinds are. (See ref. [3,4] and refs. given therein). It also
concerns choices. Choices of specific subsets of maximal chains
from a non-tree poset specifically obtained
starting from the Fibonacci rabbits` tree.\\
 Several figures illustrate the exposition of statements - the derivation of the
recurrence itself included.

\vspace{0.7cm}

\section{Introduction}
 There are various classical interpretations of
binomial coefficients,  Stirling numbers, and the $q$- Gaussian
coefficients . Recently Kovanlina  have discovered  [3,4] a
 simple and natural unified combinatorial interpretation  of all
 of them in terms of object  selection from weighted boxes with
 and without box repetition. So we are now in a position of the
 following recognition:

The classical, historically established standard interpretations might
be schematically presented for the sake of hint as follows:

SETS : Binomial coefficient $\Big({n \atop k}\Big )$
 , $\Big ({(n+k-1)\atop k}\Big ) $ denote number of subsets (without
  and with repetitions) - i.e. we are dealing with LATTICE of
  subsets.

SET PARTITIONS:  Stirling numbers of the second kind  $\Big\{ {n
\atop k}\Big\} $    denote number of partitions into exactly $k$
blocs - i.e. we are dealing with LATTICE of partitions.

PERMUTATION PARTITIONS : Stirling numbers of the first kind $\Big[
{n \atop k}\Big] $ denote number of permutations containing
exactly $k$ cycles

SPACES: $q$-Gaussian coefficient $\Big({n \atop k}\Big )_q $
denote number of $k$-dimensional subspaces in $n-th$ dimensional
space over Galois field  $GF(q)$ [5,6,7] i.e. we are dealing with
LATTICE of subspaces.

\vspace{3mm}

 \textbf{Before Konvalina combinatorial interpretation}.

\vspace{3mm}

 Algebraic similarities of the above classes of situations provided
 Rota [5] and Goldman and Rota [8,9] with an incentive to start the algebraic
unification that captures the intrinsic properties of these
numbers. The binomial coefficients, Stirling numbers and Gaussian
coefficients appear then as the coefficients in the characteristic
polynomials of geometrical lattices  [5] (see also [10] for the
subset-subspace analogy). The generalized coefficients [3]  are
called Whitney numbers of the first (characteristic polynomials)
and the second kind (rank polynomials).

\vspace{3mm}

 \textbf{Konvalina combinatorial interpretation}

 All these cases above and the case of Gaussian coefficients of the first
 kind $q^{\Big({n \atop 2}\Big )}\Big({n \atop k}\Big )_q $ are given unified
 Konvalina combinatorial interpretation in terms  of  the \textbf{generalized} binomial
 coefficients of the first and of the second kind  (see: [3,4] ).

\textbf{Unknowns ?} As for the distinguished [11,12,13,14]
Fibonomial coefficient defined below - no combinatorial
interpretation was  known till  today now to the present author .
The aim of this note is to promote a long  time waited for
 - classical in the spirit - combinatorial interpretation of Fibonomial coefficients.
 Namely we propose following [1,2] such a partial ordered set that the Fibonomial coefficients
  count the number of specific finite ``$birth-selfsimilar$''  sub-posets of an infinite
  locally finite not of binomial type , non-tree poset naturally related to
  the Fibonacci tree of rabbits growth process. This partial ordered set is defined equivalently
  via $\zeta$  characteristic matrix of partial order relation from its Hasse diagram.
  The classical scheme to be continued through  "Fibonomials" interpretation is the following:

POSET :   Fibonomial coefficient $\left( \begin{array}{c} n\\k\end{array}
\right)_{F}$  is the number of  ``birth-selfsimilar'' subposets.

\vspace{3mm}

\section{Combinatorial Interpretation}
It pays to get used to write $q$ or $\psi$ extensions of binomial symbols in
 mnemonic convenient  upside down notation [16,17] .
\begin{equation}\label{eq1}
\psi_n\equiv n_\psi , x_{\psi}\equiv \psi(x)\equiv\psi_x ,
 n_\psi!=n_\psi(n-1)_\psi!, n>0 ,
\end{equation}
\begin{equation}\label{eq2}
x_{\psi}^{\underline{k}}=x_{\psi}(x-1)_\psi(x-2)_{\psi}...(x-k+1)_{\psi}
\end{equation}
\begin{equation}\label{eq3}
x_{\psi}(x-1)_{\psi}...(x-k+1)_{\psi}= \psi(x)
\psi(x-1)...\psi(x-k-1) .
\end{equation}
You may consult [16,17] for further development and profit from
the use of this notation . So also here we use this upside down
convention for Fibonomial  coefficients:

\vspace{3mm}
$$
\left( \begin{array}{c} n\\k\end{array}
\right)_{F}=\frac{F_{n}!}{F_{k}!F_{n-k}!}\equiv
\frac{n_{F}^{\underline{k}}}{k_{F}!},\quad n_{F}\equiv F_{n}\neq 0, $$

\noindent where we make an analogy driven [16,17]  identifications $(n>0)$:
$$
n_{F}!\equiv n_{F}(n-1)_{F}(n-2)_{F}(n-3)_{F}\ldots 2_{F}1_{F};$$
$$0_{F}!=1;\quad n_{F}^{\underline{k}}=n_{F}(n-1)_{F}\ldots (n-k+1)_{F}. $$

\noindent This is the specification of the notation from [16] for the
purpose  Fibonomial Calculus case (see Example 2.1 in [17]).

\vspace{3mm}

    Let us now define the partially ordered infinite set  $P$.   We
shall label its vertices  by pairs of coordinates: ${\langle i , j
\rangle} \in {N \times N}$.   Vertices  show up in layers
("generations") of $N \times N$  along the
recurrently-subsequently emerging $s-th$ levels $\Phi_s$ where -
note! $s\in N$  i.e.

\begin{defn}
$$\Phi_s =\{\langle j, s\rangle 1\leq j \leq s_F\}, {s\in N}. $$
\end{defn}

We shall refer to $\Phi_s$  as to  the set of vertices at the
$s-th$ level. The population of the  $k-th$ level ("generation" )
counts  $k_F$  different member vertices.

\vspace{2mm}

Here down a disposal of vertices on $\Phi_k$ levels is visualized.

\vspace{5mm}

$---\Uparrow-----\Uparrow----up --Fibonacci---stairs--\star--k-th-level$\\

$---- and ----- so ---- on ---- up    --- \Uparrow ----------$\\
$\star \star \star \star \star \star \star \star \star \star \star \star \star \star \star \star \star \star\star \star \star \star \star \star \star \star \star \star \star \star \star \star \star \star \star \star \star \star \star \star \star \star \star \star --\star \star \star \star\star10-th-level$\\
$\star \star \star \star \star \star \star \star \star \star \star \star \star \star \star \star \star \star \star \star \star \star \star \star \star \star \star \star \star \star \star\star\star\star----------- 9-th-level$\\
$\star \star \star \star \star \star \star \star \star \star \star \star \star \star \star \star \star \star \star \star \star------------------8-th-level$\\
$\star \star \star \star \star \star \star \star \star \star \star \star \star -------------------------7-th-level$\\
$\star \star \star \star \star \star \star \star-----------------------------6-th-level$\\
$\star \star \star \star \star ---------------------------------5-th-level$\\
$\star \star \star ---------------------------------- 4-th-level$\\
$\star \star -----------------------------------3-rd-level  $ \\
$\star ------------------------------------ 2-nd-level$\\
$\star ----------------------------------- 1-st-level$\\

 \vspace{2mm}

   \textbf{Figure 0. The $s-th$ levels in $ N \times N $}

 \vspace{2mm}

Accompanying to the set $V$ of vertices - the set $E$ of edges
where here down ${p,q,s}\in N $ we obtain the Hasse diagram .
Namely
\begin{defn}

$$P=\langle V,E\rangle ,\\ V=\bigcup_{1\leq p}\Phi_p ,\\ E
=\{\langle\langle j , p\rangle ,\langle q ,(p+1) \rangle \rangle\}
, 1 \leq j \leq {p_F} , 1\leq q \leq {(p+1)_F}.$$
\end{defn}
\begin{defn}
The prototype (to be copied) cobweb sub-poset is : $P_m =
\bigcup_{1\leq s\leq m}\Phi_s.$
\end{defn}

In reference [2] a partially ordered infinite set $P$ was
introduced  via  descriptive picture of  its Hasse diagram. Indeed
, we may picture out the partially ordered infinite set $P$ from
the  Definition $1$ with help of the sub-poset $P_{m}$ ({\it
rooted at $F_{1}$ level of the poset}) to be continued then ad
infinitum in now obvious way as seen from the  $Fig.1$ of $P_{5}$
below. It looks like the Fibonacci rabbits` tree with a specific
``cobweb''.

\vspace{2mm}

\begin{center}

\includegraphics[width=75mm]{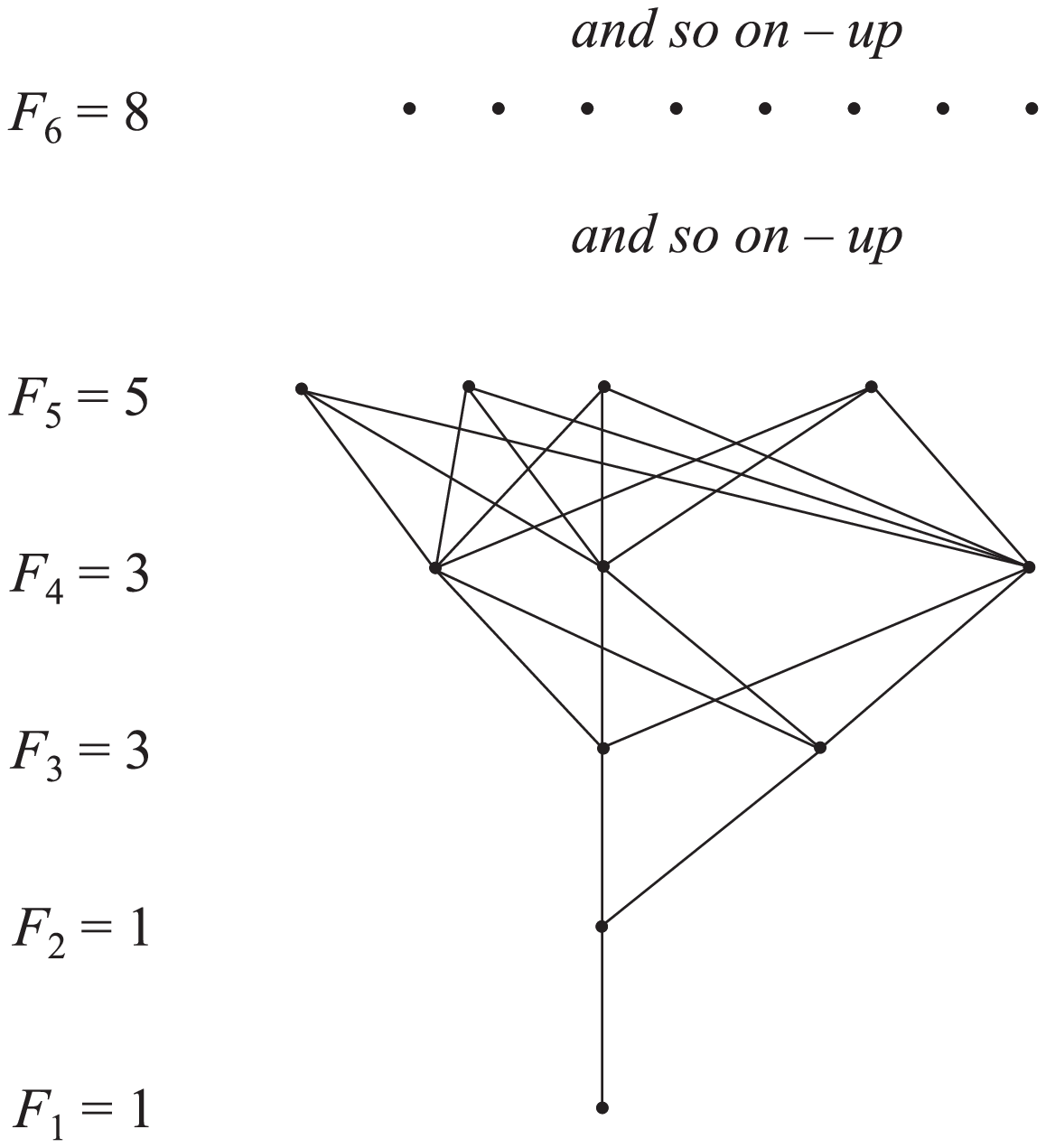}

\vspace{2mm}

\noindent {\small Fig.~1. Combinatorial interpretation of
Fibonomial coefficients.} \end{center}

\vspace{2mm} \noindent As seen above the $Fig.1$. displays the
rule of the construction of the  Fibonacci "cobweb"  poset. It is
being visualized clearly while defining this cobweb poset  $P$
with help of its incidence matrix . The incidence  $\zeta$
function [5,6,7] matrix representing uniquely just this cobweb
poset $P$ has the staircase structure correspondent with
$"cobwebed"$ Fibonacci Tree i.e. a Hasse diagram of the particular
partial order relation under consideration. This is seen below on
the Fig.$2$.

\vspace{2mm}

$$ \left[\begin{array}{cccccccccccccccc}
 1 & 1 & 1 & 1 & 1 & 1 & 1 & 1 & 1 & 1 & 1 & 1 & 1 & 1 & 1 & \cdots\\
 0 & 1 & 1 & 1 & 1 & 1 & 1 & 1 & 1 & 1 & 1 & 1 & 1 & 1 & 1 & \cdots\\
 0 & 0 & 1 & 0 & 1 & 1 & 1 & 1 & 1 & 1 & 1 & 1 & 1 & 1 & 1 & \cdots\\
 0 & 0 & 0 & 1 & 1 & 1 & 1 & 1 & 1 & 1 & 1 & 1 & 1 & 1 & 1 & \cdots\\
 0 & 0 & 0 & 0 & 1 & 0 & 0 & 1 & 1 & 1 & 1 & 1 & 1 & 1 & 1 & \cdots\\
 0 & 0 & 0 & 0 & 0 & 1 & 0 & 1 & 1 & 1 & 1 & 1 & 1 & 1 & 1 & \cdots\\
 0 & 0 & 0 & 0 & 0 & 0 & 1 & 1 & 1 & 1 & 1 & 1 & 1 & 1 & 1 & \cdots\\
 0 & 0 & 0 & 0 & 0 & 0 & 0 & 1 & 0 & 0 & 0 & 0 & 1 & 1 & 1 & \cdots\\
 0 & 0 & 0 & 0 & 0 & 0 & 0 & 0 & 1 & 0 & 0 & 0 & 1 & 1 & 1 & \cdots\\
 0 & 0 & 0 & 0 & 0 & 0 & 0 & 0 & 0 & 1 & 0 & 0 & 0 & 1 & 1 & \cdots\\
 0 & 0 & 0 & 0 & 0 & 0 & 0 & 0 & 0 & 0 & 1 & 0 & 1 & 1 & 1 & \cdots\\
 0 & 0 & 0 & 0 & 0 & 0 & 0 & 0 & 0 & 0 & 0 & 1 & 1 & 1 & 1 & \cdots\\
 0 & 0 & 0 & 0 & 0 & 0 & 0 & 0 & 0 & 0 & 0 & 0 & 1 & 0 & 0 & \cdots\\
 0 & 0 & 0 & 0 & 0 & 0 & 0 & 0 & 0 & 0 & 0 & 0 & 0 & 1 & 0 & \cdots\\
 0 & 0 & 0 & 0 & 0 & 0 & 0 & 0 & 0 & 0 & 0 & 0 & 0 & 0 & 1 & \cdots\\
 . & . & . & . & . & . & . & . & . & . & . & . & . & . & . & \cdots\\
 \end{array}\right]$$

 \vspace{2mm}

\textbf{Figure 2.  The  staircase structure  of  incidence matrix
$\zeta$ }

\vspace{2mm}

\textbf{Description} In the $k-th$ row  the "cob re-viewer"
encounters $F_k -1$ zeros right to the diagonal value $1$ thus
getting a picture of descending down to infinity led by diagonal
direction with use of growing in size cobweb Fibonacci staircase
tiled and build of the only up the diagonal zeros
 - note these are forbiddance zeros (they code no edge links along $k-th$ levels
 ("generations") of the "$cobwebed$"  Fibonacci rabbits tree from [2].

\vspace{2mm}

This staircase structure of incidence [6,7] matrix $\zeta$ which
equivalently defines uniquely  this particular cobweb poset was
being recovered right from the Definition 1 and illustrative Hasse
diagram in Fig.$1$ of Fibonacci cobweb  poset. Let us say it again
- if one decides to define the poset $P$ by incidence matrix
$\zeta$  then must arrives at $\zeta$ with this easily
recognizable staircase-like structure of zeros in the upper part
of this upper triangle incidence matrix $\zeta$ just right from
the picture (see [18,1]).

Let us recall [5,6,7] that $\zeta$ is being defined for any poset
as follows
 $(x,y \in P)$:
$$   \zeta (x,y)=\left\{ \begin{array}{cl} 1&for\ x \leq y,\\0&otherwise.
\end{array} \right. $$

\vspace{2mm}

The above $\zeta$ characteristic matrix of the partial order
relation in $P$ has been expressed explicitly in [1] in terms of
the infinite Kronecker delta  matrix $\delta$ from incidence
algebra $I(P)$ [5,6,7] as follows: ($\langle x,y\rangle \in
{N\times \Phi}$ ).

\vspace{2mm}
                   $$\zeta = \zeta_1 - \zeta_0 $$ where for $x,y \in N$,
$$\zeta_1(x,y) = \sum_{k\geq 0} \delta (x+k,y)$$  while

$$\zeta_0(x,y)=\sum_{k\geq 0} \sum_{s\geq 1}
\delta (x,F_{s+1}+k)\sum_{1\leq r\leq (F_s
-k-1)}\delta(k+F_{s+1}+r,y).$$\\
 Naturally

$$   \delta(x,y)=\left\{ \begin{array}{cl} 1&for\ x = y,\\0&otherwise.

\end{array} \right. $$

\textbf{Important.} The knowledge of $\zeta$  matrix explicit form
enables one [6,7] to construct (count) via standard algorithms
[6,7]  the M{\"{o}}bius matrix $\mu =\zeta^{-1} $ and other
typical elements of incidence algebra perfectly suitable for
calculating number of chains, of maximal chains etc. in finite
sub-posets of $P$. Right from the definition of $P$ via its Hasse
diagram  here now obvious observations follow .

\vspace{2mm}

\noindent {\bf Observation 1}

{\it The number of maximal chains starting from The Root  (level
$1_F$) to reach any point at the $n-th$ level  with $n_F$ vertices
is equal to $n_{F}!$}.

\vspace{2mm}

\noindent {\bf Observation 2} $(k>0)$

{\it The number of maximal chains rooted in any fixed vertex at
the $k-th$ level reaching the $n-th$ level with $n_F$ vertices is
equal to $n_{F}^{\underline{m}}$, where $m+k=n.$ }

\vspace{2mm}

Indeed. Denote the number of ways to get along
maximal chains from a fixed point in $\Phi_k $to$  \Rightarrow  \Phi_n , n>k$ with the symbol\\
  $$[\Phi_k \rightarrow \Phi_n]$$
  then obviously we have :\\
           $$[\Phi_1 \rightarrow \Phi_n]= n_F!$$ and
$$[\Phi_1 \rightarrow \Phi_k]\times [\Phi_k\rightarrow \Phi_n]=
[\Phi_1 \rightarrow \Phi_n].$$

\vspace{2mm}

  In order to find out the combinatorial interpretation of
  Fibonomial coefficients let us make use of the following active analogy
  of making a choice of finite sub-poset $P_m(k)_r$ max-disjoint copies rooted at
  $k-th$ level at a fixed root $\langle r,k \rangle, 1 \leq r \leq k_F $  and
  ending  at corresponding vertices (sub-cobweb leafs) up at the $n-th$
  level where the leafs live ($n=k+m$).
Explanation: max-disjoint means that sub-posets looked upon as
families of \textbf{maximal} chains are disjoint.
  In coordinate system we define the cobweb sub-poset $P_m(k)_r$ as follows:

\begin{defn}
Let   $1 \leq r \leq {k_F} , 0 \leq s \leq {m-1}.$ Let $\langle
k,r\rangle \oplus \Pi$ denotes the shift of the set $\Pi$ with the
vector $\langle k,r\rangle$ Let $\Phi_o = \{\langle 0,0
\rangle\}.$
 Then we define:\\
$P_m(k)_r=\langle V_m(k)_r,  E_m(k)_r\rangle , V_m(k)_r= \langle
k,r\rangle \oplus \bigcup_{0\leq s\leq m}\Phi_s ,\\
E_m(k)_r =\{\langle\langle {(r+j)},(k+s)\rangle,\langle (r+i)
,(k+s+1)\rangle \rangle , 1\leq (r+j) \leq (k+s)_F , 1\leq {(r+i)}
\leq {(k+s+1)_F}\}.$
\end{defn}
\textbf{Observe}    $P_m(1)_1 = \langle V_m(1)_1, E_m(1)_1 \rangle
\equiv\langle V_m, E_m\rangle \equiv P_m .$ Hence $V_m(k)_r
=\langle k,r\rangle \oplus V_m $. Here, let us recall: $P_m$ is
the sub-poset of $P$ rooted at the $1-th$ level consisting of all
intermediate level vertices up to $m-th$ level ones - those from
$\Phi_m$ included $(See: Fig.1.)$.

\vspace{2mm}

\textbf{A newly k-th level born sub-cob browsing.}\\
 Consider now the following behavior of a sub-cob useful animal
moving from any given point of the $F_k$ "generation level"  of
the poset up and then up... It behaves as it has been born right
there and can reach at first $F_2$ vertices-points up, then $F_3$
points up , $F_4$ up... and so on - thus climbing up to the level
$F_{k+m} = F_n$ of the poset  $P$.  It can see - as its Great
Ancestor at the Source Root $F_1-th$ Level- and then potentially
follow-  one of its own thus  accessible max-disjoint
\textit{copy} of sub-poset $P_m(k)_r$ .

One of many of such max-disjoint sub-posets $P_m$`s copies rooted
at a fixed point of  the $k-th$ level might be then found as a
good genetic choice  to start thus limited maximal chains
forwarding up. How many choices can be made?

\vspace{5mm}

\noindent {\bf Observation 3} $(k>0)$

{\it Let  $n = k+m$. The number of max-disjoint  sub-poset $P_{m}$
copies rooted at any fixed point at the $k-th$ level  and ending
at the n-th level  is equal to}
$$\frac{n_{F}^{\underline{m}}}{m_{F}!} =\left( \begin{array}{c} n\\m\end{array}
\right)_{F}= \left( \begin{array}{c} n\\k\end{array} \right)_{F}=
\frac{n_{F}^{\underline{k}}}{k_{F}!}. $$

\vspace{2mm}

Indeed. Consider the number of all max-disjoint sub-poset
$P_m(k)_r$ \textit{copies} rooted at the fixed vertex $\langle
(r+j),k \rangle , 1\leq (r+j) \leq k_F $. Denote this number with
the symbol

$$ \left( \begin{array}{c} n\\k\end{array}\right)_{F}$$

then obviously you have :

\begin{equation}
\frac {[\Phi_1 \rightarrow \Phi_n]}{[\Phi_1 \rightarrow \Phi_k]} =
\left( \begin{array}{c} n\\k\end{array}\right)_{F} \times [\Phi_1
\rightarrow \Phi_m]
\end{equation}

Indeed. It is enough to notice that $[\Phi_1 \rightarrow \Phi_m]$
counts the number of maximal chains in any copy of the $P_m$.
\renewcommand{\thesubsection}{\arabic{subsection}.}

\vspace{3mm}

\section{Does Konvalina like interpretation  of objects $F$- selections from weighted boxed exist?}

Binomial enumeration or finite operator calculus of Roman-Rota and Others is now the
standard tool of combinatorial analysis. The corresponding $q$-binomial calculus ($q$-calculus
- for short) is also the basis of much numerous  applications (see [19,20] for
altogether couple of thousands of respective references via enumeration and links).
In this context Konvalina unified binomial coefficients look intriguing and much promising.
The idea of  $F$-binomial or Fibonomial finite operator calculus (see Example 2.1 in [17])
consists of specification of the general scheme  - (see: [16,17] and references also to
Ward, Steffensen ,Viskov , Markowsky and others - therein)- specification  via the choice
 of the Fibonacci sequence to be sequence defining the generalized binomiality of polynomial
 bases involved (see Example 2.1 in [17]).Till now however we had been  lacking alike combinatorial
 interpretation of Fibonomial coefficients. We hope that this note would help not only via
 Observations above but also due to coming next- observation where recurrence relation
 for Finonomial coefficients is derived (recognized) thanks to its combinatorial interpretation.

\vspace{2mm}

\noindent{\bf Observation 4} $(k>0)$  ,  (combinatorial
interpretation of the recurrence)

The following known [11,14] recurrences hold
$$\left( \begin{array}{c} {n+1}\\k\end{array}
\right)_{F}=  F_{k-1} \left( \begin{array}{c} n\\k\end{array}
\right)_{F} + F_{n-k+2} \left( \begin{array}{c} n\\{k-1}\end{array}
\right)_{F}$$
or equivalently
$$\left( \begin{array}{c} {n+1}\\k\end{array}
\right)_{F}=  F_{k+1} \left( \begin{array}{c} n\\k\end{array}
\right)_{F} + F_{n-k} \left( \begin{array}{c} n\\{k-1}\end{array}
\right)_{F}$$ where
$$\left( \begin{array}{c} n\\0\end{array}
\right)_{F}= 1 , \left( \begin{array}{c} 0\\k\end{array}
\right)_{F}=0 , $$ \vspace{2mm} due to the recognition that we are
dealing with two disjoint classes in $P_{(n+1)}$ (n =k+m).
\textbf{The first one} for which
 $$ F_{k+1} \left( \begin{array}{c} n\\k\end{array}
\right)_{F}$$ equals to  $F_{k+1}$  times number of max-disjoint
copies of $P_m$ - rooted at a fixed point on the $k-th$ level (see
Interpretation below) and

 \vspace{2mm}

\textbf{the second one} for which

\begin{equation}\label{eq4}
 F_{n-k} \left( \begin{array}{c} n\\{k-1}\end{array}\right)_{F}=
 F_{n-k} \left( \begin{array}{c} n\\{n-k+1}\end{array}\right)_{F}
\end{equation}

equals to $F_{n-k}$ times number of max-disjoint $P_m`s$ copies-
rooted ata fixed point at the $(k-1)-th$ level  and ending at the
$n-th$ level - see Interpretation below.

 \vspace{2mm}

 \textbf{Interpretation}   $(k>0)$

 \vspace{2mm}

$\diamond \diamond \diamond \diamond \diamond \diamond \diamond
\diamond \diamond \diamond F_{(n-k+2)} \diamond \diamond \diamond
\diamond \star\star\star\star\star\star\star \star-- \star
F_{(n+1)}-F_{(n-k+2)} \star\star\star\star\star\star\star
$\\

$\diamond\diamond\diamond\diamond\diamond\diamond\diamond\diamond\diamond\diamond\diamond\diamond\diamond\diamond\diamond\diamond\diamond\diamond\diamond\diamond\diamond---\diamond\diamond-----\star \star \star \star----n-th--level \star --- \star\star $\\

$\diamond\diamond\diamond\diamond\diamond\diamond\diamond\diamond\diamond\diamond\diamond\diamond\diamond\diamond\diamond\diamond\diamond\diamond\diamond\diamond\diamond\star\star\star\star\star\star\star\star\star\star\star\star\star\star\star\star\star\star\star\star\star\star\star--\star\star\star\star\star\star\star--\star\star\star\star\star\star\star$\\

$\diamond \diamond \diamond \diamond \diamond \diamond \diamond \diamond \diamond \diamond \diamond \diamond \diamond \star \star \star \star \star\star \star \star \star \star \star \star \star \star \star \star \star \star \star \star \star \star \star \star \star \star \star \star \star \star \star \star \star \star \star --10-th-level$ \\

$\diamond \diamond \diamond \diamond \diamond \diamond \diamond \diamond \star \star \star \star \star \star \star \star \star \star \star \star \star \star \star \star \star \star \star \star \star \star \star\star\star\star-------- 9-th-level$\\

$\diamond \diamond \diamond \diamond \diamond \star \star \star \star \star \star \star \star \star \star \star \star \star\star\star\star\star\star----------------8-th-level$\\

$\diamond\diamond\diamond \star \star \star \star \star \star \star \star \star \star -------------------------7-th-level$\\

$\diamond\diamond \star \star \star \star \star \star-----------------------------6-th-level$\\

$\diamond \star \star \star \star ---------------------------------5-th-level$\\

$\diamond \star \star ---------------------------------- 4-th-level$\\

$\star \star -----------------------------------3-rd-level  $ \\

$\star ------------------------------------ 2-nd-level$\\

$\star ----------------------------------- 1-st-level$\\

 \vspace{2mm}

   \textbf{Figure 3. The Diamond choice - two disjoint classes in $P_{(m+1)}(k)_r$.}

 \vspace{2mm}

The Fig.3 illustrates how the two disjoint classes referred to in
Observation 4  come into existence ($r=1, k=4$). First: every
cobweb sub-poset has the "trunk" of length $\geq$ one (in the
$Fig.3$ it is the extreme left maximal chain).  From any selected
root-vertex in $k-th$ level $F_{(k+1)}$ trunks may be continued in
$F_{k+1}$ ways. A trunk of the $P_{m+1}$ copy  being chosen - for
example the set of vertices $\langle1,s \rangle , k \leq s \leq
{(n+1)}$) in the case of diamond cobweb poset selected in Fig. 3 -
the resulting sub-cobweb ends with $F_{(n-k+2)}$ diamond vertices
("leafs") at $(n+1)-th$ level. The max-disjoint copies when
shifted (in $F_{k+1}$ ways - each ) up and correspondingly
completed by - with the ultimate rightist maximal chain ending -
lacking part of now $P_{m+1}`s$ copy become max-disjoint copies
rooted at $k-th$ level and ending at the $(n+1)th$ level. This
gives

 $$ F_{k+1}\left(\begin{array}{c}
n\\k\end{array}\right)_{F} ,$$

what constitutes  the first summand of the corresponding
recurrence.

 \vspace{2mm}

Consider then the non-$\Phi_k$ level (then to be shifted) choice
of the vertex .

 \vspace{5mm}

$\diamond \diamond \diamond \diamond \diamond \diamond \diamond
\diamond \diamond \diamond F_{(n+1)}-F_{(n-k+3)} \diamond \diamond
\diamond \diamond \star\star\star\star\star\star\star --\star
\otimes F_{(n-k+3)}
\otimes\otimes\otimes---\otimes\otimes\otimes\otimes\otimes$\\

$\diamond\diamond\diamond\diamond\diamond\diamond\diamond\diamond\diamond\diamond\diamond\diamond\diamond\diamond\diamond\diamond\diamond\diamond\diamond\diamond---\diamond\diamond-----\star \star \star -- \otimes--n-th--level \otimes --- \otimes\otimes $\\

$\diamond\diamond\diamond\diamond\diamond\diamond\diamond\diamond\diamond\diamond\diamond\diamond\diamond\diamond\diamond\diamond\diamond\diamond\diamond\diamond\diamond\star\star\star\star\star\star\star\star\star\star\star\star\star\star\star\star\star--\otimes\otimes\otimes\otimes\otimes\otimes\otimes--\otimes\otimes\otimes\otimes\otimes\otimes\otimes$\\

$\diamond \diamond \diamond \diamond \diamond \diamond \diamond \diamond \diamond \diamond \diamond \diamond \diamond \star \star \star \star \star\star \star \star -- \star \star \otimes\otimes\otimes\otimes\otimes\otimes\otimes\otimes\otimes\otimes\otimes\otimes\otimes\otimes\otimes\otimes\otimes\otimes\otimes\otimes\otimes 10-th-level$ \\

$\diamond \diamond \diamond \diamond \diamond \diamond \diamond \diamond \star \star \star \star \star \star \star  \star  \star  \star \star \otimes\otimes\otimes\otimes\otimes\otimes\otimes\otimes\otimes\otimes\otimes\otimes\otimes----------- 9-th-level$\\

$\diamond \diamond \diamond \diamond \diamond \star \star  \star  \star  \star \star \star \star \otimes\otimes\otimes\otimes\otimes\otimes\otimes\otimes------------------8-th-level$\\

$\diamond\diamond\diamond \star \star \star  \star  \star \otimes\otimes\otimes\otimes\otimes -----------------------7-th-level$\\

$\diamond\diamond \star \star \star \otimes\otimes\otimes----------------------------6-th-level$\\

$\diamond \star \star \otimes\otimes ---------------------------------5-th-level$\\

$\diamond \star \otimes---------------------------------- 4-th-level$\\

$\star \otimes--------------------------------3-rd-level$\\

$\otimes-------------------------------- 2-nd-level$\\

$\star ----------------------------------- 1-st-level$\\

 \vspace{2mm}

   \textbf{Figure 4. The non-diamond choice - two disjoint classes in $P_{(n+1)}$.}

 \vspace{2mm}

\vspace{2mm}

The Fig.4 continues to illustrate how the two disjoint classes
referred to in Observation 4  are introduced. Now - what we do we
choose a vertex-root $\otimes$ in $(k-1)-th$ level in one of
$F_{(k-1)}$ ways. A trunk being chosen - say of the $\otimes$
cobweb sub-poset in Fig. 4 - it ends with $F_{(n-k+2)}$ $\otimes$
vertices ("leafs") at $n-th$ level. The number of $\otimes$
max-disjoint copies of $P_{m+1}$ rooted at a fixed point at
$(k-1)-th level$ is equal to  $ \left(\begin{array}{c}
n\\k-1\end{array}\right)_{F}$. These max-disjoint copies  become
max-disjoint copies rooted at $k-th$ level and ending at the
$(n+1)-th$ level when shifted up the $k-th$ and rooted at the same
"diamond" root of the first choice and then  correspondingly
completed by - with the ultimate leftist maximal chain ending -
lacking part of  $P_{m+1}`s$ max-disjoint copies with all other
copies rooted there at $k-th$ level and ending at the $(n+1)-th$
level. The number of thus obtained max-disjoint copies is equal to

$$F_{n-k} \left( \begin{array}{c} n\\{k-1}\end{array}
\right)_{F}.$$

 All together this gives the number
of all max-disjoint cobweb sub-posets ending at $\Phi_{(n+1)}$
while starting from a fixed point of  $\Phi_k$ level. It is equal
to the sum of cases in the two disjoint classes  i.e.

$$\left( \begin{array}{c} {n+1}\\k\end{array}
\right)_{F}=  F_{k+1} \left( \begin{array}{c} n\\k\end{array}
\right)_{F} + F_{n-k} \left( \begin{array}{c} n\\{k-1}\end{array}
\right)_{F}.$$

\vspace{5mm}

 In this connection the intriguing  question arises : May
one extend-apply somehow Konvalina theorem [3,4] below so as to
encompass also Fibonomial case under investigation ?

In [3,4] Konvalina considers  $n$ distinct boxes labeled with $i\in {[n]} , [n] \equiv \{1,...n\}$
such that each of $i-th$ box contains  $w_i$ distinct objects. John Konvalina uses the convention
$1\leq w_1 \leq w_2\leq ...\leq w_n$. Vector $N^n \ni \vec w =(w_1, w_2 ,...,w_n)$ is the weigh
 vector then. Along with Konvalina  considerations we have from [3]:

\vspace{3mm} \textbf{The Konvalina Theorem 1}

\vspace{2mm}

Let  $\vec w =(w_1, w_2 ,...,w_n)$ where $ 1\leq w_1 \leq w_2 \leq ...\leq w_n$.  Then

\textbf{I.} $$  C_{k}^n(\vec w) = C_{k}^{n-1}(\vec w) + w_n C_{k-1}^{n-1}(\vec w)$$

\textbf{II.}$$  S_{k}^n(\vec w) = S_{k}^{n-1}(\vec w) + w_n S_{k-1}^{n-1}(\vec w).$$

Here  $C_{k}^n(\vec w)$  denotes the generalized binomial coefficient of the first
kind with weight $\vec w$ and  it is the number of ways to select $k$ objects from
$k$ (necessarily distinct !) of the $n$ boxes with constrains as follows :  choose
$k$ distinct labeled boxes $$ i_1<i_2 <...<i_k $$  and then choose one object
 from each of the $k$ distinct boxes selected.  Naturally  one then has [3]

$$C_{k}^n(\vec w) = \sum_{1\leq i_1< i_2 < ...< i_k\leq n} w_{i_1}w_{i_2}...w_{i_k} . $$

Complementarily   $S_{k}^n(\vec w)$ denotes the generalized binomial coefficient of the second
kind with weight $\vec w$ and  it is the number of ways to select $k$ objects from
$k$ (not necessarily  distinct)  of the $n$ boxes with constrains as follows [3]:  choose $k$
not necessarily distinct labeled boxes $$ i_1\leq i_2 \leq ...\leq i_k$$  and then choose
one object from  each of the $k$ (not necessarily  distinct) boxes selected.  Obviously
one then has

$$S_{k}^n(\vec w) = \sum_{1\leq i_1\leq i_2 \leq ...\leq i_k\leq n} w_{i_1}w_{i_2}...w_{i_k}.$$

here the natural question arises : how are we to extend  Konvalina
theorem [3,4] above so as to encompass also Fibonomial case under
investigation ?

\vspace{2mm}

\textbf{Information I }: \textit about the preprint [21] entitled
\textit{Determinants, Paths,  and Plane Partitions} by Ira M.
Gessel,  X. G. Viennot [21].

Right after  Theorem $25$ - Section 10 , page 24 in [21])-
relating  the number $ N(R)$ of nonintersecting $k$-paths to
Fibonomial coefficients via $q$-weighted type counting formula-
the authors express their  wish worthy to be quoted: "\textit{it
would be nice to have a more natural interpretation then the one
we have given}"... " \textit{ R. Stanley has asked if there is a
binomial poset associated with the Fibonomial coefficients}..." -
Well. The cobweb locally finite infinite poset  by Kwasniewski
from [15,18,1,2] is not of binomial type. Recent incidence algebra
origin arguments [22] seem to make us not to expect  binomial type
poset come into the game.
    The $q$-weighted type counting formula from [21]gives rise to
an interesting definition of Fibonomial coefficients all together
with its interpretation in terms of nonintersecting $k$-paths due
to the properties of binomial determinant. Namely ,following [21]
let us consider points  $P_i= \langle\ 0,-i \rangle$ and $Q_i =
\langle\-n+i,-n+i \rangle$. Let  $R=\{r_1 < r_2 < ...<r_k\}\equiv
R(\vec r)$ be a subset of $\{0,1,...,n\}\equiv[n+1]$. Let $N(R)$
denotes the number of non-intersecting $k$-paths from $\langle
P_{r_1},...,P_{r_k}\rangle$ to $\langle
Q_{r_1},...,Q_{r_k}\rangle$ . Then $det \left( \begin{array}{c}
{r_i}\\n-r_{k+1-j}\end{array} \right) $ = $N(R)$.  The
$q$-weighted type counting formula from [21]then for $q=1$ means
that

$$\left( \begin{array}{c} {n+1}\\k\end{array}
\right)_{F}= \sum_{R(\vec r)}N(R) .$$

    In view of [21] another question arises - what is the relation like between
these two: Gessel and Viennot [21] non-intersecting $k$-paths and
cobweb sub-poset [18,2,3] points of view?

\vspace{2mm}

\textbf{Information II }: \textit on the partial ordered poset and
Fibonacci numbers paper [23] by Istvan Beck. The author of [23]
shows that  $F_n$ equals to the number of of ideals in a simple
poset called "fence" . This allows Him to infer via combinatorial
reasoning the identities :

$$F(n) = F(k) F(n + 1 - k) + F(k - 1)F(n -k)$$
$$F(n)= F(k-1) F(n + 1 -( k-1)) + F(k - 2) F(n -(k-1)).$$
 A straightforward application
of these above is the confirmation - just by checking - the
intriguing validity of recurrence relation for Fibonomial
coefficients . As we perhaps might learn from this note coming to
the end - both the Fibonomial coefficients as well as their
recurrence relation  are interpretable along the classical
historically established manner referring to the number of
objects` choices - this time these are partially ordered sub-sets
here called the cobweb sub-posets - the effect of the diligent
spider`s spinning of the maximal chains cobweb   during the
arduous day spent on the infinite Fibonacci rabbits` growth tree.

 \vspace{2mm}

\textbf{{\bf} Historical Memoir Remark}
The Fibonacci sequence origin is attributed and referred to the first
edition (lost)  of  ``Liber abaci'' (1202) by Leonardo Fibonacci  [Pisano]
(see second edition from 1228 reproduced as Il Liber Abaci di Leonardo
Pisano publicato secondo la lezione Codice Maglibeciano by  Baldassarre
Boncompagni  in Scritti di Leonardo Pisano  vol. 1, (1857) Rome).

\vspace{2mm}

\textbf{{\bf} Historical Quotation Remark } As accurately noticed
by Knuth and Wilf in [14]  the recurrent relations for Fibonomial
coefficients appeared already in $1878$  Lukas work [11]. In our
opinion - Lucas`s Th\'eorie des fonctions num\'eriques simplement
p\'eriodiques is the far more non-accidental context for binomial
and binomial-type coefficients - Fibonomial coefficients included.

While studying this mentioned important and inspiring paper by
Knuth and Wilf [14] and in the connection with a context of this
note a question raised by the authors with respect to their
formula (15) is worthy to be repeated : \textit{Is there a
"natural" interpretation....} - May be then fences from [23]  or
cobweb posets or  ... "Natural" naturally might have many
effective faces ...

\vspace{2mm}

\begin
{thebibliography}{99}
\parskip 0pt

\bibitem{1}
A. K. Kwa\'sniewski, {\it More on Combinatorial interpretation of
the Fibonomial coefficients} ArXiv: math.CO/0402344 v1 23 Feb 2004

\bibitem{2}
A. K. Kwasniewski {\it Information on combinatorial interpretation
of Fibonomial coefficients }   Bull. Soc. Sci. Lett. Lodz Ser.
Rech. Deform. 53, Ser. Rech.Deform. {\bf 42} (2003): 39-41 ArXiv:
math.CO/0402291   v1 18 Feb 2004

\bibitem{3}
J. Konvalina , {\it Generalized binomial coefficients and the
subset-subspace  problem } , Adv. in Appl. Math. {\bf 21}  (1998)
: 228-240

\bibitem{4}
J. Konvalina , {\it A Unified Interpretation of the Binomial
Coefficients, the Stirling Numbers and the Gaussian Coefficients}
The American Mathematical Monthly {\bf 107}(2000):901-910

\bibitem{5}
 Gian-Carlo Rota "On the Foundations of Combinatorial Theory,
I. Theory  of  Möbius Functions";   Z. Wahrscheinlichkeitstheorie
und Verw. Gebiete, vol.2 ,  1964 , pp.340-368.

\bibitem{6}
 E. Spiegel, Ch. J. O`Donnell  {\it Incidence algebras}  Marcel
Dekker, Inc. Basel $1997$ .

\bibitem{7}
Richard P. Stanley {\it Enumerative Combinatorics}  {\bf I},
Wadsworth and Brooks Cole Advanced Books  and Software, Monterey
California, $1986$

\bibitem{8}
 Goldman J.  Rota G-C. {\it The Number of Subspaces in a vector
space} in  Recent Progress in Combinatorics (W. Tutte, Ed.):
75-83, Academic Press, San Diego, $1969$, see: ( J. Kung, Ed.)
"Gian Carlo Rota on Combinatorics"  Birkhäuser, Boston
(1995):217-225

\bibitem{9}
Goldman J.  Rota G-C. {\it On the Foundations of Combinatorial
Theory IV; finite-vector spaces and Eulerian generating functions}
Studies in Appl. Math. {\bf 49} (1970): 239-258

\bibitem{10}
J. Kung   {\it The subset-subspace analogy} ( J. Kung, Ed.)
"Gian Carlo Rota on Combinatorics"  Birkhäuser, Boston (1995):277-283

\bibitem{11}
E. Lucas, {\it Th\'eorie des fonctions num\'eriques simplement
p\'eriodiques}, American Journal of Mathematics {\bf 1} (1878):
184--240; (Translated from the French by Sidney Kravitz), Ed. D.
Lind, Fibonacci Association, 1969.

\bibitem{12}
 G. Fonten\'e, {\it G\'en\'eralisation d`une formule connue},
Nouvelles Annales de Math\'ematiques (4) {\bf 15} (1915), 112.

\bibitem{13}
H. W. Gould,   {\it The bracket function and Fonten\'e-Ward
generalized binomial coefficients with applications to Fibonomial
coefficients}, The Fibonacci Quarterly {\bf 7} (1969), 23--40.

\bibitem{14}
D. E. Knuth, H. S. Wilf  {\it The Power of a Prime that Divides a
Generalized Binomial Coefficient} J. Reine Angev. Math. {\bf 396}
(1989): 212-219

\bibitem{15}
A. K. Kwa\'sniewski, {\it More on Combinatorial interpretation of
Fibonomial coefficients},   Inst. Comp. Sci.  UwB/Preprint no. 56,
November 2003.

\bibitem{16}
A. K. Kwa\'sniewski,   {\it Towards  $\psi$-extension of finite
operator calculus of Rota}, Rep. Math. Phys. {\bf 47} no. 4
(2001), 305--342.    ArXiv: math.CO/0402078  2004

\bibitem{17}
A. K. Kwa\'sniewski, {\it On simple characterizations of Sheffer
$\Psi$-polynomials and related propositions of the calculus of sequences},
Bull.  Soc.  Sci.  Lettres  \L \'od\'z {\bf 52},S\'er. Rech. D\'eform.
 {\bf 36} (2002), 45--65. ArXiv: math.CO/0312397  $2003$

\bibitem{18}
A. K. Kwa\'sniewski,{\it Combinatorial interpretation of Fibonomial coefficients}
 Inst. Comp. Sci.  UwB Preprint No  {\bf 52} (November $2003$ )

\bibitem{19}
T. Ernst  , {\it The History of q-Calculus and a new Method },
 http://www.math.uu.se/~thomas/Lics.pdf 19 December (2001),
(Licentiate Thesis). U. U. D. M. Report (2000).

\bibitem{20}
A. K. Kwa\'sniewski,{\it First Contact Remarks on Umbra Difference Calculus References Streams}
 Inst. Comp. Sci.  UwB Preprint No  {\bf 63} (January $2004$ )

\bibitem{21}
Ira M. Gessel,  X. G. Viennot{\it Determinant Paths and Plane
Partitions  } preprint  (1992)
http://citeseer.nj.nec.com/gessel89determinants.html

\bibitem{22}
A.K.Kwasniewski {\it The second part of on duality triads`
paper-On fibonomial and other triangles  versus  duality triads}
 Bull. Soc. Sci. Lett. Lodz Ser. Rech. Deform. 53, Ser. Rech.
Deform. {\bf 42} (2003): 27 -37  ArXiv: math.GM/0402288 v1 18 Feb.
$2004$

\bibitem{23}
I. Beck  {\it Partial Orders and the Fibonacci Numbers} The
Fibonacci Quarterly {\bf 26} (1990): 172-174 .

\end{thebibliography}



\end{document}